\documentclass[12pt]{article}
\usepackage{amsfonts,amssymb,amsmath}
\usepackage{epsfig}
\usepackage[all]{xy}
\usepackage[makeroom]{cancel}

\newcommand\CC{\mathbb C}

\renewcommand{\Re}{\mathop{\mathrm{Re}}}
\renewcommand{\Im}{\mathop{\mathrm{Im}}}

\newcommand\beq{\begin{equation}}
\newcommand\eeq{\end{equation}}

\newtheorem{theorem}{Theorem}

\newtheorem{remark}{Remark}

\begin{document}

\title{{Moutard transform for the generalized analytic functions}
\thanks{The main part of the work was fulfilled during the visit of the 
first author to the Centre de Math\'ematiques Appliqu\'ees of \'Ecole Polytechnique
in October 2015. The first author was partially supported by the Russian Foundation for Basic Research, 
grant 13-01-12469 ofi-m2, by the program ``Leading scientific schools'' (grant NSh-4833.2014.1), 
by the program ``Fundamental problems of nonlinear dynamics''. }}

\author{P.G. Grinevich
\thanks{L.D. Landau Institute for Theoretical Physics,
pr. Akademika Semenova 1a, 
Chernogolovka, 142432, Russia; Lomonosov Moscow State University, Faculty of Mechanics and Mathematics, 
Russia, 119991, Moscow, GSP-1, Leninskiye Gory 1, Main Building;
Moscow Institute of Physics and Technology, 9 Institutskiy per., Dolgoprudny,
Moscow Region, 141700, Russia; e-mail: pgg@landau.ac.ru} \and R.G. Novikov\thanks
{CNRS (UMR 7641), Centre de Math\'ematiques Appliqu\'ees, 
\'Ecole Polytechnique, 91128, Palaiseau, France; IEPT RAS, 117997, Moscow, Russia;
e-mail: novikov@cmap.polytechnique.fr}}
\date{}
\maketitle
\begin{abstract}
We construct a Moutard-type transform for the generalized analytic functions. The first theorems 
and the first explicit examples in this connection are given.
\end{abstract}

\section{Introduction}
We consider the equation
\begin{equation}
\label{eq:ga1}
\partial_{\bar z} \psi = u \bar \psi \ \ \mbox{in} \ \ D\subseteq\CC,
\end{equation}
where $u=u(z)$ is a given function in $D$, $D$ is open in $\CC$.
The functions $\psi=\psi(z)$ satisfying (\ref{eq:ga1}) are known as generalized 
analytic functions in $D$. In this article the notation $f=f(z)$ does not mean that 
$f$ is holomorphic. In the literature it is usually assumed that
\begin{align}
\label{eq:cond1}
&u \in L_p(D),& \  &p>2,& \  &\mbox{if} \ \ D \ \ \mbox{is bounded},& \hspace{4cm}\\
\label{eq:cond2}
&u \in L_{p,2}(\CC),&  \ &p>2,&  \ &\mbox{if} \ \ D=\CC,& \hspace{4cm}
\end{align}
where  
\begin{align}
&L_{p,\nu}(\CC)\ \ \mbox{denotes complex-valued functions} \ \ u \ \ \mbox{such that}
\nonumber \\
&u\in  L_{p}(D_1), \ \ u_{\nu}\in  L_{p}(D_1), \ \ \mbox{where} \ \ u_{\nu}(z)=
\frac{1}{|z|^{\nu}} u\left(\frac{1}{z}\right),\\
&D_1=\{z\in\CC\, :\, |z|\le 1 \}. \nonumber
\end{align}
The theory of generalized analytic functions  is presented in \cite{Bers}, \cite{Vek}.
However, to the best of our knowledge, algebraic Moutrad-type transforms, going back 
to \cite{Mout} were not yet considered in the framework of this theory. 
On the other hand, the Moutard-type transforms are successfully used in  
studies of integrable models in dimension 2+1 and in the spectral theory in 
dimension 2; see \cite{DGNS}, \cite{MatvSal}, \cite{NS}, \cite{NTT}-\cite{TT} 
and references therein. 

The Moutard-type transforms correspond in quadratures to coefficients and related 
solutions of an appropriate linear PDE (and its conjugate) on the plane  new coefficients 
and related solutions for this PDE  (and its conjugate). In this article we show that this 
approach is applicable to equation (\ref{eq:ga1}). In particular, it allows to study equation  
(\ref{eq:ga1}) in some important cases when $u$ has 
strong singularities (e.g. contour poles), and, as a corollary, assumptions  (\ref{eq:cond1}),
(\ref{eq:cond2}) are not valid at all. It is quite likely that the known methods of the 
generalized analytic functions theory dot not admit appropriate generalizations for the 
aforementioned cases.

In addition, our interest to the generalized analytic functions with strong singularities was, 
in particular, strongly motivated by studies going back to \cite{GN}.

Note that the present work was stimulated by recent articles \cite{Taim1}, \cite{Taim2} 
by I.A. Taimanov on Moutard-type transforms for two-dimensional 
Dirac operators with applications to integrable systems in dimension 2+1 and to differential 
geometry.

The results of the present work can be summarized as follows: 

\begin{itemize}
\item We construct Moutard-type transforms for equation (\ref{eq:ga1}) considering this 
equation as a particular case of the two-dimensional Dirac system
\begin{equation}
\label{eq:dirac1}
\partial_{\bar z} \psi_1= u \psi_2, \ \ \partial_{z} \psi_2= v \psi_1 \ \ \mbox{in} \ \ D\subseteq\CC,
\end{equation}
where $D$ is the same that in (\ref{eq:ga1}). Related algebraic and analytic results in this connection 
are presented in Section~\ref{sec:2}.  

\item Explicit examples of generalized analytic functions obtained by the Moutard-type 
transforms from usual holomorphic functions are given in Section~\ref{sec:3}. These 
examples include generalized analytic functions with contour poles.
\end{itemize}

\section{Moutard transform}
\label{sec:2}
We consider the two-dimensional Dirac system (\ref{eq:dirac1}) and the conjugate system 
\begin{equation}
\label{eq:dirac2}
\partial_{\bar z} \psi^+_1= -v \psi^+_2, \ \ \partial_{z} \psi^+_2= -u \psi^+_1 \ \ \mbox{in} \ \ D\subseteq\CC.
\end{equation} 

If 
\begin{equation}
\label{eq:reality1}
v=\bar u,
\end{equation}
 then system (\ref{eq:dirac1}) is reduces to (\ref{eq:ga1}). More precisely, in this case:
\begin{align}
&\mbox{if} \  \psi \  \mbox{satisfies} \  (\ref{eq:ga1}), \ \mbox{then} 
\left[\begin{array}{c} \psi_1 \\ \psi_2 \end{array} \right] = 
\left[\begin{array}{c} \psi \\ \bar\psi \end{array} \right] 
\  \mbox{satisfies} \ (\ref{eq:dirac1});
\label{eq:reality2} \\
&\mbox{if} \ \left[\begin{array}{c} \psi_1 \\ \psi_2 \end{array} \right] 
\  \mbox{satisfies} \ (\ref{eq:dirac1}), \ \mbox{then} \ \psi=\frac{1}{2}(\psi_1+\bar\psi_2) 
\ \mbox{and} \  \psi=\frac{1}{2i}(\psi_1-\bar\psi_2) \  \  \mbox{satisfy} \  (\ref{eq:ga1}).
\nonumber
\end{align}

In addition, if (\ref{eq:reality1}) holds, then the conjugate system (\ref{eq:dirac2}) is reduced to the equation 
\begin{equation}
\label{eq:ga2}
\partial_{\bar z} \psi^+ = -\bar u \bar \psi^+ \ \ \mbox{in} \ \ D\subseteq\CC.
\end{equation}

In the theory of generalized analytic functions equation (\ref{eq:ga2}) is known as the conjugate 
equation to the initial equation (\ref{eq:ga1}); see \cite{Vek}. 

Let 
\begin{equation}
\label{eq:dirac:sol1}
\vec\psi(j)=\left[\begin{array}{c} \psi_1(z,j) \\ \psi_2(z,j) \end{array} \right],
\ \ 
\vec\psi^+(j)=\left[\begin{array}{c} \psi^+_1(z,j) \\ \psi^+_2(z,j) \end{array} \right], \ \ j=1,\ldots,N,
\end{equation}
denote a set of fixed solutions of  systems (\ref{eq:dirac1}) and  (\ref{eq:dirac2}), respectively.

In addition, for systems  (\ref{eq:dirac1}) and  (\ref{eq:dirac2}), we consider the potentials 
$\omega_{j,k}=\omega_{j,k}(z)$ defined as follows:
\begin{equation}
\label{eq:dirac3} 
\partial_{z} \omega_{j,k} =\psi_1(j)\psi^+_1(k), \ \ 
\partial_{\bar z} \omega_{j,k} =-\psi_2(j)\psi^+_2(k) \ \ \mbox{in} \ \ 
D\subseteq\CC, \ \ j,k=1,\ldots,N,
\end{equation}
under the assumption that $D$ is simply connected. Note that 
$$
\partial_{\bar z} \partial_{z}\omega_{j,k} = 
(\partial_{\bar z} \psi_1(j))\ \psi^+_1(k) +
\psi_1(j)\ (\partial_{\bar z} \psi^+_1(k))
=u \psi_2(j)\ \psi^+_1(k) - v \psi_1(j)\ \psi^+_2(k),
$$
$$
\partial_{z} \partial_{\bar z} \omega_{j,k} = 
-(\partial_{z} \psi_2(j))\ \psi^+_2(k) -
\psi_2(j)\ (\partial_{z} \psi^+_2(k))=
-v  \psi_1(j)\ \psi^+_2(k) +u
\psi_2(j)\ \psi^+_1(k),
$$
and, thus,
\begin{equation}
\label{eq:dirac4} 
\partial_{\bar z} \partial_{z}\omega_{j,k} = 
\partial_{z} \partial_{\bar z} \omega_{j,k}.
\end{equation}
Therefore, definitions (\ref{eq:dirac3}) are self-consistent, and the 
integration constants may depend on the concrete situation.

Let us introduce $N\times N$ matrix $\Omega=(\Omega_{k,j}(z))=(\omega_{j,k}(z)) $:
\begin{equation}
\label{eq:omega1}
\Omega=\begin{bmatrix} 
\omega_{1,1} & \omega_{2,1} & \ldots & \omega_{N,1}  \\
\omega_{1,2} & \omega_{2,2} & \ldots & \omega_{N,2}  \\
\vdots & \vdots & \ddots & \vdots \\
\omega_{1,N} & \omega_{2,N} & \ldots & \omega_{N,N}  \\
\end{bmatrix}.
\end{equation}
\begin{theorem} \label{th:1}
Suppose that $\vec\psi(j)$, $\vec\psi^+(j)$, $j=1,\ldots,N$, 
(see (\ref{eq:dirac:sol1})) are 
formal solutions of systems (\ref{eq:dirac1}), (\ref{eq:dirac2}), respectively, 
for given coefficients $u=u(z)$, $v=v(z)$, and that $\det\Omega\ne0$, where $\Omega$
is defined by (\ref{eq:dirac3})-(\ref{eq:omega1}). Let the transform  
\begin{align}
&\{u,v\}\rightarrow \{\tilde u,\tilde v\},  \label{eq:moutard1} \\ 
&\{\vec\psi(0),\vec\psi^+(0)\} \rightarrow \{\vec{\tilde\psi}(0), 
\vec{\tilde\psi}^+(0) \} \label{eq:moutard2}
\end{align}
be defined as follows:
\begin{equation}
\label{eq:moutard3}
\tilde u = u + \begin{bmatrix}\psi_1(1) & \ldots & \psi_1(N) \end{bmatrix} 
 \Omega^{-1} \begin{bmatrix} \psi^+_2(1) \\ \vdots \\ \psi^+_2(N)  \end{bmatrix},
\end{equation}
\begin{equation}
\label{eq:moutard4}
\tilde v= v - \begin{bmatrix}\psi_2(1) & \ldots & \psi_2(N) \end{bmatrix} 
 \Omega^{-1}  \begin{bmatrix} \psi^+_1(1) \\ \vdots \\ \psi^+_1(N)  \end{bmatrix} ,
\end{equation}
\begin{equation}
\label{eq:moutard5}
\begin{bmatrix}\tilde\psi_1(0) \\ \tilde\psi_2(0) \end{bmatrix}=
\begin{bmatrix}\psi_1(0) \\ \psi_2(0) \end{bmatrix}-
\begin{bmatrix}\psi_1(1) & \ldots & \psi_1(N) \\ 
\psi_2(1) & \ldots & \psi_2(N) \end{bmatrix} 
\Omega^{-1}
\begin{bmatrix}\omega_{0,1} \\ \vdots \\ \omega_{0,N} 
\end{bmatrix},
\end{equation}
\begin{equation}
\label{eq:moutard6}
\begin{bmatrix}\tilde\psi^+_1(0) \\ \tilde\psi^+_2(0) \end{bmatrix}=
\begin{bmatrix}\psi^+_1(0) \\ \psi^+_2(0) \end{bmatrix}-
\begin{bmatrix}\psi^+_1(1) & \ldots & \psi^+_1(N) \\ 
\psi^+_2(1) & \ldots & \psi^+_2(N) \end{bmatrix} 
(\Omega^{-1})^t 
\begin{bmatrix}\omega_{1,0} \\ \vdots \\ \omega_{N,0} 
\end{bmatrix},
\end{equation}
where $\vec\psi(0)$,  $\vec\psi^+(0)$ are some formal solutions to systems 
(\ref{eq:dirac1}), (\ref{eq:dirac2}), respectively, $\omega_{0,j}$ and  
$\omega_{j,0}$ are defined as in (\ref{eq:dirac3}), 
but in terms of $\vec\psi(0)$,
$\vec\psi^+(j)$ and  $\vec\psi(j)$, $\vec\psi^+(0)$, respectively, and the 
symbol $t$  in (\ref{eq:moutard6}) stands for the matrix transposition.
Then the transformed functions $\vec{\tilde\psi}(0)$,  $\vec{\tilde\psi}^+(0)$ 
solve the transformed Dirac equations:
\begin{align}
\label{eq:dirac5}
&\partial_{\bar z} \tilde\psi_1(0)= \tilde u \tilde\psi_2(0),& \ \ 
&\partial_{z} \tilde\psi_2(0)= \tilde v \tilde\psi_1(0)& \ \ &\mbox{in} \ \ D,\\
\label{eq:dirac6}
&\partial_{\bar z} \tilde\psi^+_1(0)= -\tilde v \tilde\psi^+_2(0),& \ \ 
&\partial_{z} \tilde\psi^+_2(0)= -\tilde u \tilde\psi^+_1(0)& \ \ &\mbox{in} \ \ D.
\end{align}
\end{theorem}

Historically, transformations like (\ref{eq:moutard1})-(\ref{eq:moutard2}) 
go back to the paper \cite{Mout}, and are known in the literature as Moutard 
transforms or Moutard-type transforms. For the two-dimensional Dirac operators 
the Moutard transforms were considered many times in the literature, see, for 
example, \cite{MatvSal}, \cite{Taim1}, \cite{Taim2} and references therein. Nevertheless, 
for completeness of exposition, we present the proof of Theorem~\ref{th:1} in 
Section~\ref{sec:4}. 

The point is that the Moutard-type transform (\ref{eq:moutard1}),(\ref{eq:moutard2}) admits the following 
reduction to the case of equations  (\ref{eq:ga1}),  (\ref{eq:ga2}), i.e. to the case of 
generalized analytic functions.

Let 
\begin{equation}
\label{eq:ga:sol1}
\psi(j)=\psi(z,j), \ \ \psi^+(j)=\psi^+(z,j), \ \ j=1,\ldots,N,
\end{equation}
denote a set of fixed solutions of  (\ref{eq:ga1}) and (\ref{eq:ga2}), respectively. Then we can
define $\omega_{j,k}=\omega_{j,k}(z)$ as imaginary-valued functions satisfying
\begin{equation}
\label{eq:ga:kern1}
\partial_{z} \omega_{j,k} =\psi(j)\psi^+(k), \ \ 
\partial_{\bar z} \omega_{j,k} =-\overline{\psi(j)\psi^+(k)} \ \ \mbox{in} \ \ 
D, \ \ j,k=1,\ldots,N,
\end{equation}
where (\ref{eq:ga:kern1}) is the reduction of (\ref{eq:dirac3}) corresponding to  
\begin{equation}
\label{eq:ga:kern2}
\psi_1(j)=\psi(j), \ \ \psi_2(j)=\overline{\psi(j)}, \ \
\psi_1^+(j)=\psi^+(j), \ \ \psi^+_2(j)=\overline{\psi^+(j)}.
\end{equation}

\begin{remark}
It is interesting to note that real-valued function $\phi=\frac{\omega_{j,k}}{2i}$, $i=\sqrt{-1}$, 
for some fixed $j$, $k$, is known in the generalized analytic functions theory as a potential for equation 
(\ref{eq:ga1}); see \cite{Vek}.
\end{remark}

\begin{theorem} \label{th:2}
Suppose that $\psi(j)$, $\psi^+(j)$, $j=1,\ldots,N$, 
(see (\ref{eq:ga:sol1})) are 
formal solutions of equations (\ref{eq:ga1}), (\ref{eq:ga2}), respectively, 
for a given coefficient $u=u(z)$, and that $\det\Omega\ne0$, where $\Omega$
is defined according  to (\ref{eq:omega1}), (\ref{eq:ga:kern1}) with 
imaginary-valued $\omega_{j,k}$. Let the transform  
\begin{align}
&u \rightarrow \tilde u,  \label{eq:ga:moutard1} \\ 
&\{\psi(0),\psi^+(0)\} \rightarrow \{{\tilde\psi}(0),{\tilde\psi}^+(0) \}  
\label{eq:ga:moutard2}
\end{align}
be defined as follows:
\begin{equation}
\label{eq:ga:moutard3}
\tilde u = u + \begin{bmatrix}\psi(1) & \ldots & \psi(N) \end{bmatrix} 
 \Omega^{-1} \begin{bmatrix} \overline{\psi^+(1)} \\ \vdots \\ 
\overline{\psi^+(N)}  \end{bmatrix},
\end{equation}
\begin{equation}
\label{eq:ga:moutard5}
\tilde\psi(0) =
\psi(0) -
\begin{bmatrix}\psi(1) & \ldots & \psi(N) \end{bmatrix} 
\Omega^{-1}
\begin{bmatrix}\omega_{0,1} \\ \vdots \\ \omega_{0,N} 
\end{bmatrix},
\end{equation}
\begin{equation}
\label{eq:ga:moutard6}
\tilde\psi^+(0)= \psi^+(0) -
\begin{bmatrix}\psi^+(1) & \ldots & \psi^+(N) \end{bmatrix} 
(\Omega^{-1})^t 
\begin{bmatrix}\omega_{1,0} \\ \vdots \\ \omega_{N,0} 
\end{bmatrix},
\end{equation}
where $\psi(0)$,  $\psi^+(0)$ are some formal solutions to equations 
(\ref{eq:ga1}), (\ref{eq:ga2}), respectively, $\omega_{0,j}$ and  
$\omega_{j,0}$ are imaginary-valued and are defined as in (\ref{eq:ga:kern1}), 
but in terms of $\psi(0)$, $\psi^+(j)$ and  $\psi(j)$, $\psi^+(0)$, 
respectively, and $t$ in (\ref{eq:ga:moutard6}) stands for the matrix 
transposition.
Then the transformed functions ${\tilde\psi}(0)$,  ${\tilde\psi}^+(0)$ 
solve the transformed generalized-analytic function equations:
\begin{align}
\label{eq:ga:moutard7}
&\partial_{\bar z} \tilde\psi(0)= \tilde u\, \overline{\psi(0)}& \ \ &\mbox{in} \ \ D,\\
\label{eq:ga:moutard8}
&\partial_{\bar z} \tilde\psi^+(0)= -\overline{\tilde u}\, 
\overline{\tilde\psi^+(0)}& \ \ 
&\mbox{in} \ \ D.
\end{align}
\end{theorem}

To our best knowledge, Moutard-type type transforms like 
(\ref{eq:ga:moutard1}), (\ref{eq:ga:moutard2}) were not yet considered in the 
generalized-analytic function theory. 

Theorem~\ref{th:2} follows from Theorem~\ref{th:1} in the framework of 
reductions (\ref{eq:reality1}), (\ref{eq:ga:kern2}). In this case using
also that $\overline{\Omega}=-\Omega$ one can see that the transform 
(\ref{eq:moutard1}), (\ref{eq:moutard2}) preserves symmetries  
(\ref{eq:reality1}), (\ref{eq:ga:kern2}).

The algebraic result of Theorem~\ref{th:2} admits, in particular, the following analytic realization:

\begin{theorem}
\label{th:3}  
Let $D$ be an open simply-connected bounded domain in $\CC$ with $C^1$-boundary and $u$ satisfy 
(\ref{eq:cond1}). Let $\psi(j)\in W^{1,p}(D)$ and $\psi^+(j)\in W^{1,p}(D)$, $j=1,\ldots,N$, be solutions
of equations (\ref{eq:ga1}) and (\ref{eq:ga2}), respectively, and $\det\Omega\ne0$ in $D\cup \partial D$, 
where $\Omega$ is defined according  to (\ref{eq:omega1}), (\ref{eq:ga:kern1}) with 
imaginary-valued $\omega_{j,k}$. Let the Moutard transform be defined according the formulas 
(\ref{eq:ga:moutard1})-(\ref{eq:ga:moutard6}) as in Theorem~\ref{th:2}. 

Then the transformed coefficient $\tilde u$ satisfied (\ref{eq:cond1}) as well as the initial $u$.
\end{theorem}

In Theorem~\ref{th:3}, $W^{1,p}$ denotes the standard Sobolev space. 

Theorem~\ref{th:3} follows from formula (\ref{eq:ga:moutard3}) and from the properties that 
$\psi(j)$, $\psi^+(j)$, $j=1,\ldots,N$, and $\omega_{j,k}$,  $1\le j,k\le N$, and $\det\Omega$ are 
continuous in $D\cup \partial D$. Note that, in order to obtain the aforementioned continuity of 
$\psi(j)$, $\psi^+(j)$, we use that, under our assumptions on $D$, if $f\in W^{1,p}(D)$ then 
$f\in C^{\alpha}(D\cup\partial D)$, where $C^{\alpha}$ denotes the standard H\"older space, 
$\alpha=(p-2)/p$ (see \cite{AF}). Finally, the continuity of 
$\omega_{j,k}$ follows from the continuity of $\psi(j)$, $\psi^+(k)$ and from formula (\ref{eq:dirac3}).

\section{Explicit  examples}
\label{sec:3}

In this section we present explicit examples illustrating the Moutard-type transforms for 
generalized analytic functions, i.e. transforms (\ref{eq:ga:moutard1})-(\ref{eq:ga:moutard6}).

Let 
\begin{align}
\label{eq:ex:1}
&u\equiv 0, \ \ N=1, \ \ \psi(1)=f(z), \ \  \psi^+(1)=f^+(z), \\
&\psi(0)=\psi(z,0), \ \ \psi^+(0)=\psi^+(z,0), \nonumber
\end{align}
where $f$, $f^+$ and $\psi(0)$, $\psi^+(0)$  are holomorphic on $D$. Then:
\begin{enumerate}
\item For imaginary-valued $\omega_{1,1}$ and  $\omega_{0,1}$,  $\omega_{1,0}$  
of (\ref{eq:ga:kern1}) and  (\ref{eq:ga:moutard5}), (\ref{eq:ga:moutard6}), we have
\begin{align}
\label{eq:ex:2}
&\omega_{1,1}(z)=F_{1,1}(z)-\overline{F_{1,1}(z)}+c_{1,1},\\
&\omega_{0,1}(z)=F_{0,1}(z)-\overline{F_{0,1}(z)}+c_{0,1}, \ \ 
\omega_{1,0}(z)=F_{1,0}(z)-\overline{F_{1,0}(z)}+c_{1,0},
\nonumber
\end{align}
where $F_{1,1}$ and $F_{0,1}$, $F_{1,0}$ are holomorphic functions on $D$ such that 
\begin{align}
\label{eq:ex:3}
&\partial_z F_{1,1}(z) =f(z)f^+(z),\\
&\partial_z F_{0,1}(z) =\psi(z,0)f^+(z),\ \ \partial_z F_{1,0}(z) =f(z)\psi^+(z,0),
\nonumber
\end{align}
and $c_{1,1}$, $c_{0,1}$, $c_{1,0}$  are pure imaginary constants;
\item Formulas (\ref{eq:ga:moutard3})-(\ref{eq:ga:moutard6}) take the form:
\begin{equation}
\label{eq:ex:4}
\tilde u(z) = \frac{f(z) \overline{f^+(z)}}{\omega_{1,1}(z)},
\end{equation}
\begin{equation}
\label{eq:ex:5}
\tilde \psi(z,0) = \psi(z,0)-f(z)\,   \frac{\omega_{0,1}(z)}{\omega_{1,1}(z)},
\end{equation}
\begin{equation}
\label{eq:ex:6}
\tilde \psi^+(z,0) = \psi^+(z,0)-f^+(z)\,  \frac{\omega_{1,0}(z)}{\omega_{1,1}(z)}.
\end{equation}
\end{enumerate}
In addition, equations (\ref{eq:ga:moutard7}), (\ref{eq:ga:moutard8}) for $\tilde\psi(0)$,
$\tilde\psi^+(0)$ are nontrivial already for 
\begin{equation}
\label{eq:ex:6.1}
f(z)=C\ne0, \ \ f^+(z)=C^+\ne0,
\end{equation}
where $C$, $C^+$ are complex constants. In particular, in this case
\begin{equation}
\label{eq:ex:7}
\tilde u(z) = \frac{C \overline{C^+}}{2i \Im(CC^+ z)+ c_{1,1}}.
\end{equation}

One can see that this $\tilde u$ satisfies (\ref{eq:cond1}) if 
$(D\cup\partial D) \subset \Lambda^{+}$ or 
$(D\cup\partial D)\subset \Lambda^{-}$, where
\begin{align}
\label{eq:ex:8}
&\Lambda^{+} = \{ z\in\CC: \Im\left(2 CC^+ z+ c_{1,1}\right) >0 \} , \\ 
&\Lambda^{-} = \{ z\in\CC: \Im\left(2 CC^+ z+ c_{1,1}\right) <0 \}  .\nonumber 
\end{align}
However, this $\tilde u$ does not satisfy (\ref{eq:cond1}) at all if $D\cap\mathcal{L}\ne\emptyset$,
where 
\begin{equation}
\label{eq:ex:9}
\mathcal{L} = \{ z\in\CC: \Im\left(2 CC^+ z+ c_{1,1}\right) =0 \}. \\ 
\end{equation}

The point is  that already the Moutard transform (\ref{eq:ga:moutard1})-(\ref{eq:ga:moutard6}) 
for the case (\ref{eq:ex:1}), (\ref{eq:ex:6.1}) yields generalized analytic functions with contour
poles: with poles on $\mathcal{L}$ in this particular case. 

Next, we consider the example when  
\begin{align}
\label{eq:ex:10}
&u\equiv 0, \ \ N=2, \ \ \psi(1)=f_1, \ \  \psi^+(1)=f_1^+, \ \   
\psi(2)=f_2, \ \  \psi^+(2)=f_2^+,   \\
&\psi(0)=\psi(z,0), \ \ \psi^+(0)=\psi^+(z,0), \nonumber
\end{align}
where $f_1$, $f_1^+$, $f_2$, $f_2^+$ are complex constants, and $\psi(0)$,  $\psi^+(0)$ are 
holomorphic functions on $D$. In this case: 
\begin{enumerate}
\item For imaginary-valued $\omega_{j,k}$ of (\ref{eq:ga:kern1}),  
(\ref{eq:ga:moutard5}), (\ref{eq:ga:moutard6}), $j,k\in\{0,1,2\}$, we have
\begin{align}
\label{eq:ex:12}
&\omega_{j,k}(z)=f_jf^+_kz - \overline{f_jf^+_kz}+c_{j,k}, \\
&\omega_{0,k}(z)=\Psi(z)f^+_k- \overline{\Psi(z)f^+_k} +c_{0,k}, \ \ 
\omega_{j,0}(z)=f_j\Psi^+(z) -\overline{f_j\Psi^+(z)}+c_{j,0},
\nonumber
\end{align}
where $j,k=1,2$, $\Psi$, $\Psi^+$ are holomorphic functions on $D$ such that 
\begin{equation}
\label{eq:ex:13}
\partial_z \Psi(z) = \psi(z,0), \ \  \partial_z \Psi^+(z) = \psi^+(z,0),
\end{equation}
and $c_{j,k}$ are pure imaginary constants, $j,k=0,1,2$;
\item Formulas (\ref{eq:ga:moutard3})-(\ref{eq:ga:moutard6}) take the form:
\begin{equation}
\label{eq:ex:14}
\tilde u(z) = \frac{ f_{1} \overline{f^+_{1}} c_{2,2}-
f_{2} \overline{f^+_{1}} c_{1,2}-
f_{1} \overline{f^+_{2}} c_{2,1}+ f_{2} \overline{f^+_{2}}  c_{1,1}}
{\det\Omega},
\end{equation}
\begin{align}
\label{eq:ex:15}
&\tilde \psi(z,0) = \psi(z,0)- \left(\det\Omega\right)^{-1}\times  \\
&\times \left(\left[
\left( -f_{1}\overline{f_{2} f^+_{2}}+f_{2}
\overline{f_{1} f^+_{2}} \right) \overline{z}+f_{1} c_{2,2}
-f_{2}c_{1,2}\right] \omega_{0,1}(z)+\right.\nonumber\\
&\left.+\left[\left( f_{1} \overline{f_{2} f^+_{1}} 
-f_{2} \overline{f_{1} f^+_{1}}  \right) \overline{z}-f_{1}c_{2,1}
+f_{2}c_{1,1} \right] \omega_{0,2}(z) \right), \nonumber
\end{align}
\begin{align}
\label{eq:ex:16}
&\tilde \psi^+(z,0) = \psi^+(z,0)- \left(\det\Omega\right)^{-1}\times \\
&\times\left(\left[ \left( -f^+_{1}\overline{f_{2} f^+_{2}}+ f^+_{2}
\overline{f_{2} f^+_{1}} \right) \overline{z}+f^+_{1} c_{2,2}
-f^+_{2}c_{2,1}\right] \omega_{1,0}(z)+\right. \nonumber\\
&\left.+\left[\left( f^+_{1}\overline{f_{1} f^+_{2}}-f^+_{2}
\overline{f_{1} f^+_{1}} \right) \overline{z}-f^+_{1}c_{1,2}
+f^+_{2}c_{1,1} \right] \omega_{2,0}(z) \right), \nonumber
\end{align}
where
\begin{align}
\label{eq:ex:17}
&\det\Omega =2\Re\left[ \left( f_{2}f^+_{1}\overline{f_{1} f^+_{2}}-f_{1}f^+_{{1}}
\overline{f_{2} f^+_{2}} \right)z\bar z \right.+ \\
&\left.\vphantom{\overline{f_{1} f^+_{2}}-f_{1}f^+_{{1}}} 
+\left(c_{2,2} f_{1}f^+_{1} + c_{1,1}f_{2}f^+_{2}- c_{1,2} f_{2}f^+_{1} -
c_{2,1} f_{1}f^+_{2} \right) z \right] + c_{1,1}c_{2,2}-c_{2,1}c_{1,2},   \nonumber
\end{align}
\end{enumerate}
$\Omega$ is the matrix defined according to (\ref{eq:omega1}) for $N=2$.

Note that if 
\begin{align}
\label{eq:ex:18}
&c_{2,2} f_{1}f^+_{1} + c_{1,1}f_{2}f^+_{2}- c_{1,2} f_{2}f^+_{1} -
c_{2,1} f_{1}f^+_{2}=0, \\
&f_{1} \overline{f^+_{1}} c_{2,2}- f_{2} \overline{f^+_{1}} c_{1,2}-
f_{1} \overline{f^+_{2}} c_{2,1}+ f_{2} \overline{f^+_{2}}  c_{1,1} \ne 0,\nonumber\\
&\Re\left(  f_{2}f^+_{1}\overline{f_{1} f^+_{2}}-f_{1}f^+_{{1}}
\overline{f_{2} f^+_{2}} \right)\ne 0, 
\nonumber
\end{align}
then $\tilde u$ in (\ref{eq:ex:14}) is spherically-symmetric and non-trivial.
In addition, assuming that
\begin{equation}
\label{eq:ex:19}
\sigma\ge 0, \ \ \mbox{where} \ \ \sigma=-\frac{c_{1,1}c_{2,2}-c_{2,1}c_{1,2}}
{2\Re \left(  f_{2}f^+_{1}\overline{f_{1} f^+_{2}}-f_{1}f^+_{{1}}
\overline{f_{2} f^+_{2}} \right)},
\end{equation}
we have that: $\tilde u$ has a pole on 
\begin{equation}
\label{eq:ex:20}
S_r=\{z\in\CC : |z|^2=r \}, \ \ \mbox{where} \ \  r^2=\sigma;
\end{equation}
\begin{equation}
\label{eq:ex:21}
\tilde u(z)=O(|z|^{-2}) \ \ \mbox{as} \ \ |z|\rightarrow\infty.
\end{equation}

For example, assumptions (\ref{eq:ex:18}) are fulfilled if 
\begin{equation}
\label{eq:ex:22}
f_1=1, \ \ f^+_1=1, \ \ f_2=i, \ \ f^+_2=i, \ \ c_{1,1}=c_{2,2}, \ \ c_{1,2}=-c_{2,1}, \ \ 
|c_{1,1}|^2+ |c_{1,2}|^2\ne 0.
\end{equation}

One can see that the Moutard transform (\ref{eq:ga:moutard1})-(\ref{eq:ga:moutard6}) 
for the case (\ref{eq:ex:10}), (\ref{eq:ex:18}), (\ref{eq:ex:19}) yields again 
generalized analytic functions with contour poles: with poles on $S_r$ in this particular case.

Actually, our Moutard-type transforms (\ref{eq:ga:moutard1})-(\ref{eq:ga:moutard6})   
give a method for developing a proper theory of generalized analytic functions with contour 
and point poles of some natural class. This issue will be 
developed in the subsequent work. 

\section{Proof of Theorem~\ref{th:1}}
\label{sec:4}
Using (\ref{eq:dirac3}),(\ref{eq:omega1}) we obtain:
\begin{equation}
\label{eq:pth1:1}
\partial_{z}\Omega=
\begin{bmatrix} \psi^+_1(1) \\ \vdots \\ \psi^+_1(N)  \end{bmatrix}  
\begin{bmatrix} \psi_1(1) & \ldots & \psi_1(N)   \end{bmatrix},
\end{equation}
\begin{equation}
\label{eq:pth1:2}
\partial_{\bar z}\Omega=-\begin{bmatrix} \psi^+_2(1) \\ \vdots \\ \psi^+_2(N)  \end{bmatrix}  
\begin{bmatrix} \psi_2(1)  & \ldots & \psi_2(N)   \end{bmatrix},
\end{equation}
\begin{equation}
\label{eq:pth1:3}
\partial_{z}\Omega^{-1}=-\Omega^{-1}
\begin{bmatrix} \psi^+_1(1)  \\ \vdots \\ \psi^+_1(N)  \end{bmatrix}  
\begin{bmatrix} \psi_1(1) & \ldots & \psi_1(N)   \end{bmatrix}
\Omega^{-1},
\end{equation}
\begin{equation}
\label{eq:pth1:4}
\partial_{\bar z}\Omega^{-1}= \Omega^{-1} \begin{bmatrix} \psi^+_2(1) \\ \vdots \\ \psi^+_2(N)  \end{bmatrix}  
\begin{bmatrix} \psi_2(1) & \ldots & \psi_2(N)   \end{bmatrix}\Omega^{-1}.
\end{equation}

Using  (\ref{eq:moutard5}) and (\ref{eq:dirac1}), (\ref{eq:dirac3}), 
(\ref{eq:pth1:3}), (\ref{eq:pth1:4}),
 we obtain:
$$
\partial_{\bar z}\tilde\psi_1(0)= \partial_{\bar z}\psi_1(0)-
\begin{bmatrix}\partial_{\bar z}\psi_1(1) & \ldots & 
\partial_{\bar z}\psi_1(N) \end{bmatrix} 
 \Omega^{-1} 
\begin{bmatrix}\omega_{0,1} \\ \vdots \\ \omega_{0,N} 
\end{bmatrix}+
$$
$$
+\begin{bmatrix}\psi_1(1) & \ldots & \psi_1(N) \end{bmatrix} 
 \Omega^{-1} 
\begin{bmatrix}\psi^+_2(1) \\ \vdots \\ \psi^+_2(N)
\end{bmatrix}\psi_2(0)-
$$
$$
-\begin{bmatrix}\psi_1(1) & \ldots & \psi_1(N) \end{bmatrix} 
 \Omega^{-1} \begin{bmatrix} \psi^+_2(1) \\ \vdots \\ \psi^+_2(N)  \end{bmatrix}  
\begin{bmatrix} \psi_2(1) & \ldots & \psi_2(N)   \end{bmatrix}\Omega^{-1}
 \begin{bmatrix}\omega_{0,1} \\ \vdots \\ \omega_{0,N} 
\end{bmatrix}=
$$
$$
=\left(u + \begin{bmatrix}\psi_1(1) & \ldots & \psi_1(N) \end{bmatrix} 
 \Omega^{-1} \begin{bmatrix} \psi^+_2(1) \\ \vdots \\ \psi^+_2(N)  \end{bmatrix} \right)\times
$$
$$
\times\left(\psi_2(0)-\begin{bmatrix} \psi_2(1) & \ldots & \psi_2(N)   \end{bmatrix}\Omega^{-1}
 \begin{bmatrix}\omega_{0,1} \\ \vdots \\ \omega_{0,N} 
\end{bmatrix}\right)=
$$
\begin{equation}
\label{eq:pth1:5}
=\left(u + \begin{bmatrix}\psi_1(1) & \ldots & \psi_1(N) \end{bmatrix} 
 \Omega^{-1} \begin{bmatrix} \psi^+_2(1) \\ \vdots \\ \psi^+_2(N)  \end{bmatrix} \right)
\tilde\psi_2(0);
\end{equation}
$$
\partial_{z}\tilde\psi_2(0)= \partial_{z}\psi_2(0)-
\begin{bmatrix}\partial_{z}\psi_2(1) & \ldots & 
\partial_{z}\psi_2(N) \end{bmatrix} 
 \Omega^{-1} 
\begin{bmatrix}\omega_{0,1} \\ \vdots \\ \omega_{0,N} 
\end{bmatrix}-
$$
$$
-\begin{bmatrix}\psi_2(1) & \ldots & \psi_2(N) \end{bmatrix} 
 \Omega^{-1} 
\begin{bmatrix}\psi^+_1(1) \\ \vdots \\ \psi^+_1(N) \end{bmatrix}\psi_1(0)+
$$
$$
+\begin{bmatrix}\psi_2(1) & \ldots & \psi_2(N) \end{bmatrix} 
 \Omega^{-1} \begin{bmatrix} \psi^+_1(1) \\ \vdots \\ \psi^+_1(N)  \end{bmatrix}  
\begin{bmatrix} \psi_1(1) & \ldots & \psi_1(N)   \end{bmatrix}\Omega^{-1}
 \begin{bmatrix}\omega_{0,1} \\ \vdots \\ \omega_{0,N} 
\end{bmatrix}=
$$
$$
=\left(v - \begin{bmatrix}\psi_2(1) & \ldots & \psi_2(N) \end{bmatrix} 
 \Omega^{-1} \begin{bmatrix} \psi^+_1(1) \\ \vdots \\ \psi^+_1(N)  \end{bmatrix} \right)\times
$$
$$
\times\left(\psi_1(0)-\begin{bmatrix} \psi_1(1) & \ldots & \psi_1(N)   \end{bmatrix}\Omega^{-1}
 \begin{bmatrix}\omega_{0,1} \\ \vdots \\ \omega_{0,N} 
\end{bmatrix}\right)=
$$
\begin{equation}
\label{eq:pth1:6}
=\left(v - \begin{bmatrix}\psi_2(1) & \ldots & \psi_2(N) \end{bmatrix} 
 \Omega^{-1} \begin{bmatrix} \psi^+_1(1) \\ \vdots \\ \psi^+_1(N)  \end{bmatrix} \right)
\tilde\psi_1(0).
\end{equation}
Formulas (\ref{eq:pth1:5}),  (\ref{eq:pth1:6}) and   (\ref{eq:moutard3}),  (\ref{eq:moutard4})  
imply   (\ref{eq:dirac5}).

Using  (\ref{eq:moutard6}) and (\ref{eq:dirac1}), (\ref{eq:dirac3}), 
(\ref{eq:pth1:3}), (\ref{eq:pth1:4}), and the formulas
$$
\partial_{z}(\Omega^{-1})^t=(\partial_{z}\Omega^{-1})^t, \ \ 
\partial_{\bar z}(\Omega^{-1})^t=(\partial_{\bar z}\Omega^{-1})^t,
$$
 we obtain:
$$
\partial_{\bar z}\tilde\psi^+_1(0)= \partial_{\bar z}\psi^+_1(0)-
\begin{bmatrix}\partial_{\bar z}\psi^+_1(1) & \ldots & \partial_{\bar z}\psi^+_1(N) \end{bmatrix} 
 (\Omega^{-1})^t 
\begin{bmatrix}\omega_{1,0} \\ \vdots \\ \omega_{N,0} 
\end{bmatrix}+
$$
$$
+\begin{bmatrix}\psi^+_1(1) & \ldots & \psi^+_1(N) \end{bmatrix} 
 (\Omega^{-1})^t 
\begin{bmatrix}\psi_2(1) \\ \vdots \\ \psi_2(N)
\end{bmatrix}\psi^+_2(0)-
$$
$$
-\begin{bmatrix}\psi^+_1(1) & \ldots & \psi^+_1(N) \end{bmatrix} 
 (\Omega^{-1})^t \begin{bmatrix} \psi_2(1) \\ \vdots \\ \psi_2(N)  \end{bmatrix}  
\begin{bmatrix} \psi^+_2(1) & \ldots & \psi^+_2(N)   \end{bmatrix}(\Omega^{-1})^t
 \begin{bmatrix}\omega_{1,0} \\ \vdots \\ \omega_{N,0} 
\end{bmatrix}=
$$
$$
=\left(-v + \begin{bmatrix}\psi^+_1(1) & \ldots & \psi^+_1(N) \end{bmatrix} 
 (\Omega^{-1})^t \begin{bmatrix} \psi_2(1) \\ \vdots \\ \psi_2(N)  \end{bmatrix} \right)\times
$$
$$
\times\left(\psi^+_2(0)-\begin{bmatrix} \psi^+_2(1) & \ldots & \psi^+_2(N)   \end{bmatrix}(\Omega^{-1})^t
 \begin{bmatrix}\omega_{1,0} \\ \vdots \\ \omega_{N,0} 
\end{bmatrix}\right)=
$$
\begin{equation}
\label{eq:pth1:7}
=\left(-v + \begin{bmatrix}\psi^+_1(1) & \ldots & \psi^+_1(N) \end{bmatrix} 
 (\Omega^{-1})^t \begin{bmatrix} \psi_2(1) \\ \vdots \\ \psi_2(N)  \end{bmatrix} \right)
\tilde\psi^+_2(0);
\end{equation}
and
$$
\partial_{z}\tilde\psi^+_2(\lambda)= \partial_{z}\psi^+_2(\lambda)-
\begin{bmatrix}\partial_{z}\psi^+_2(1) & \ldots & 
\partial_{z}\psi^+_2(N) \end{bmatrix} 
 (\Omega^{-1})^t 
\begin{bmatrix}\omega_{1,0} \\ \vdots \\ \omega_{N,0} 
\end{bmatrix}-
$$
$$
-\begin{bmatrix}\psi^+_2(1) & \ldots & \psi^+_2(N) \end{bmatrix} 
 (\Omega^{-1})^t 
\begin{bmatrix}\psi_1(1) \\ \vdots \\ \psi_1(N) \end{bmatrix}\psi^+_1(0)+
$$
$$
+\begin{bmatrix}\psi^+_2(1) & \ldots & \psi^+_2(N) \end{bmatrix} 
 (\Omega^{-1})^t \begin{bmatrix} \psi_1(1) \\ \vdots \\ \psi_1(N)  \end{bmatrix}  
\begin{bmatrix} \psi^+_1(1) & \ldots & \psi^+_1(N)   \end{bmatrix}(\Omega^{-1})^t
 \begin{bmatrix}\omega_{1,0} \\ \vdots \\ \omega_{N,0} 
\end{bmatrix}=
$$
$$
=\left(-u - \begin{bmatrix}\psi^+_2(1) & \ldots & \psi^+_2(N) \end{bmatrix} 
 (\Omega^{-1})^t \begin{bmatrix} \psi_1(1) \\ \vdots \\ \psi_1(N)  \end{bmatrix} \right)\times
$$
$$
\times\left(\psi^+_1(0)-\begin{bmatrix} \psi^+_1(1) & \ldots & \psi^+_1(N) 
\end{bmatrix}(\Omega^{-1})^t
 \begin{bmatrix}\omega_{1,0} \\ \vdots \\ \omega_{N,0} 
\end{bmatrix}\right)=
$$
\begin{equation}
\label{eq:pth1:8}
=\left(-u - \begin{bmatrix}\psi^+_2(1) & \ldots & \psi^+_2(N) \end{bmatrix} 
 (\Omega^{-1})^t \begin{bmatrix} \psi_1(1) \\ \vdots \\ \psi_1(N)  \end{bmatrix} \right)
\tilde\psi^+_1(0).
\end{equation}
Formulas (\ref{eq:pth1:7}),  (\ref{eq:pth1:8}) and   (\ref{eq:moutard3}),  (\ref{eq:moutard4})  
imply   (\ref{eq:dirac6}). 

This completes the proof of Theorem~\ref{th:1}.

\end{document}